\documentclass{amsart}
\usepackage{amsmath}
\usepackage{amssymb}
\newtheorem{thm}{Theorem}[section]
\newtheorem{cor}[thm]{Corollary}
\newtheorem{conj}[thm]{Conjecture}
\newtheorem{lem}[thm]{Lemma}
\newtheorem{prop}[thm]{Proposition}
\newtheorem{cons}[thm]{Construction}
\theoremstyle{definition}
\newtheorem{defn}[thm]{Definition}
\theoremstyle{remark}
\newtheorem{rem}[thm]{Remark}
\newtheorem{ex}[thm]{Example}
\newtheorem{exs}[thm]{Examples}
\long\def\Thm#1{\begin{thm} #1 \end{thm}}
\long\def\Cor#1{\begin{cor} #1 \end{cor}}

\long\def\Lem#1{\begin{lem} #1 \end{lem}}
\long\def\Prop#1{\begin{prop} #1 \end{prop}}

\long\def\Rem#1{\begin{rem} #1 \end{rem}}

\def\Sect{\section}

\long\def\Ref#1#2#3#4#5#6{
\bibitem{#1}
{\rm #2,}
\textit{#3.}
{\rm #4}
\textbf{#5}
{\rm #6.}
}
\long\def\Refb#1#2#3#4{
\bibitem{#1}
{\rm #2,}
\textit{#3.}
#4.
}
\def\Zz{{\mathbb Z}}%blackboard bold Z
\def\Rr{{\mathbb R}}%blackboard bold R
\def\Ff{{\mathbb F}}

\def\phi{\varphi}
\def\into{\hookrightarrow}
\def\leq{\leqslant}
\def\geq{\geqslant}

\def\O{{\rm O}}

\def\PO{{\rm PO}}

\def\comp{\circ}
\def\semidirect{\rtimes}
\def\Sing{{\rm Sing}}
\def\Zero{{\rm Zero}}
\begin{document}

\title{On regular maps and parallel lines}

\author{M.~C.~Crabb}
\address{%
Institute of Mathematics\\
University of Aberdeen \\
Aberdeen AB24 3UE \\
UK}
\email{m.crabb@abdn.ac.uk}
\date{July 2022}
\begin{abstract}
Let $f: \Rr^{m+1}\to \Rr^{m+2^r}$, where
$2^{r-1}\leq m+1 <2^r$, be a continuous map.
Improving a recent result of Frick and Harrison, we
show that there are $4$ points $x_0,\, x_1,\, y_0,\, y_1$ in $\Rr^m$,
which are distinct if $m+1\not=2^{r-1}$, and satisfy
$x_0\not=x_1$, $y_0\not=y_1$,
$\{ x_0, x_1\} \not=\{ y_0,y_1\}$ if $m+1=2^{r-1}$,
such that the vectors $f(x_1)-f(x_0)$ and $f(y_1)-f(y_0)$
are parallel.
\end{abstract}
 %fibrewise topology
\maketitle
\Sect{Introduction}
We improve a recent result \cite[Theorems 1 and 3]{FH} of Frick and 
Harrison
using methods (perhaps more conceptual than those in \cite{FH}) 
from \cite{tver, mod2}.

Two vectors in a vector space will be called {\it parallel} if they
lie in a $1$-dimensional vector subspace;
this allows one or both vectors to be zero.
\Thm{\label{parallel}
Suppose that $m$ is a non-negative integer and let
$r\geq 1$ be the integer such that $2^{r-1}\leq m+1<2^r$. 
Let $f: \Rr^{m+1}\to \Rr^{m+2^r}$ be a continuous map.
Then

\smallskip

\par\noindent {\rm (a)\ }
if $m+1\not=2^{r-1}$,
there exist $4$ distinct points $x_0,\, x_1,\, y_0,\, y_1\in
\Rr^{m+1}$ such that the vectors $f(x_1)-f(x_0)$ and $f(y_1)-f(y_0)$
are parallel;

\smallskip

\par\noindent {\rm (b)\ }
if $m+1=2^{r-1}$,
there exist $4$ points $x_0,\, x_1,\, y_0,\, y_1\in
\Rr^{m+1}$ such that $x_0\not= x_1$, $y_0\not=y_1$,
$\{ x_0,x_1\} \not=\{ y_0,y_1\}$, and the vectors $f(x_1)-f(x_0)$ 
and $f(y_1)-f(y_0)$ are parallel.
}
At the extremes,
when $m+1=2^{r-1}$ we are considering a map 
$\Rr^{2^{r-1}} \to \Rr^{2^r+2^{r-1}-1}=\Rr^{3\cdot 2^{r-1}-1}$ 
-- the first few cases with $r=2,\, 3,\, 4,\, 5$
were established in \cite[Corollaries 2 and 4]{FH} --
and when $m+1=2^r-1$ a map
$\Rr^{2^r-1}\to \Rr^{2^{r+1}-1}=\Rr^{2(2^r-1)+1}$.

The assertion is trivial if $f$ is not injective.
For, if there are distinct points $x_0$ and $x_1$ such that
$f(x_0)=f(x_1)$, we can choose any two points $y_0$ and $y_1$
in the complement of $\{ x_0,x_1\}$, and then $f(x_1)-f(x_0)=0$
is parallel to $f(y_1)-f(y_0)$.

As a corollary we have the following special case of 
\cite[Theorem 6.16]{mod2}.
\Cor{\label{cor}
Let $f: \Rr^{m+1}\to \Rr^{m+2^r}$ be a continuous map,
where $2^{r-1}\leq m+1<2^r$.

Then there exist $4$ distinct points $x_0,\, x_1,\, x_2,\, x_3\in
\Rr^{m+1}$ such that the vectors $f(x_0)$, $f(x_1)$, $f(x_2)$,
$f(x_3)$ lie in some $2$-dimensional affine subspace of $\Rr^{m+2^r}$.
}
For, if say $x_0=y_0$ in (b),
the points $f(x_0),\, f(x_1)$ and $f(y_1)$ lie on an affine line
and we can then take $x_2=y_1$ and choose an arbitrary point
$x_3$ in the complement of $\{ x_0, x_1, x_2\}$.
\Sect{\label{prelude} Prelude}
The method will be illustrated by proving the following
classical result. (See, for example, \cite{handel}.)
\Prop{Suppose that $m\geq 0$ is a non-negative integer and that
$f: \Rr^{m+1}\to\Rr^{m+1}$ is a continuous map. 
Then there exist distinct points $x,\, y\in\Rr^{m+1}$ such that
the vectors $f(x)$ and $f(y)$ are parallel.
}
\Rem{The map $f: \Rr^{m+1}\to \Rr\oplus\Rr^{m+1}$,
$x\mapsto (1,x)$ has the property that 
$f(x),\, f(y)\in\Rr\oplus\Rr^{m+1}$
are parallel only if $x=y$.
}
Writing $V=\Rr^{m+1}$ and $W=\Rr^{n+1}$ where $n\geq 1$, 
we shall consider a continuous map $f : V\to W$.
The unit sphere in the Euclidean space $V$ is denoted by
$S(V)$.
We denote the real projective space on $V$ by $P(V)=S(V)/\{ \pm 1\}$,
and write $H_V$ (or simply $H$ if the meaning is clear) for the Hopf line bundle.
\begin{proof}
The assertion is trivially true if $m=0$.
So assume that $m\geq 1$.
We have a map
$$
\phi : X=S(V)\to V\oplus V,\quad  x\mapsto (x,-x).
$$
So $\phi (x)$ determines $2$ distinct points $\pm x$ of $V$.
The quotient of $X$ by the free action of the group
$\{ \pm 1\}$ is the real projective space $Y=P(V)$.
The mod $2$ Euler class $t=e(\lambda )$ of the Hopf
line bundle $\lambda=H_V$ over $P(V)$ determines the cohomology
ring
$$
H^*(Y;\,\Ff_2)=\Ff_2[t]/(t^{m+1}).
$$

Now the map $f$ determines a $\{ \pm 1\}$-equivariant map
$$
S(V) \to W\oplus W,\quad x\mapsto (f(x)+f(-x),f(x)-f(-x)),
$$
with $\{ \pm 1\}$ acting trivially on the first factor 
and as $\pm 1$ on the second,
and so a section of $(\Rr\oplus \lambda )\otimes W$ over
$Y$.
Lifting to $Y\times P(W)$ and projecting from the trivial bundle
$W$ to $W/H_W$, we get a section $s$ of the $2n$-dimensional 
real vector
bundle $(\Rr\oplus\lambda )\otimes (W/H_W)$ over $Y\times P(W)$.
A zero of $s$ gives a point $x\in S(V)$ and a line $L\in P(W)$
such that $f(x)+f(-x),\, f(x)-f(-x)\in L$, or,
equivalently, $f(x),\, f(-x)\in L$.

Now, as we recollect in the Appendix, the direct image homomorphism
(evaluation on the fundamental class of $P(V)$)
$$
\pi_! : H^{2n}(Y\times P(W);\,\Ff_2) \to H^n(Y;\,\Ff_2)
$$
maps the Euler class $e((\Rr\oplus\lambda )\otimes (W/H_W))$
to $w_n(-(\Rr\oplus\lambda ))$.

This Stiefel-Whitney class is the degree $n$ term of 
$(1+t)^{-1}=\sum_{0\leq i\leq m} t^i$, 
and so is non-zero if $n\leq m$.
Hence the Euler class of $(\Rr\oplus\lambda )\otimes (W/H_W)$
is zero, and the section $s$ has a zero.
\end{proof}
\Rem{A section of $(\Rr\oplus\lambda )\otimes W$ can be thought of
as a vector bundle homomorphism $\Rr\oplus \lambda \to P(V)\times W$
to the trivial bundle over $P(V)$. If this bundle homomorphism is
injective in each fibre, then $\Rr\oplus\lambda$ is included as
a subbundle of $P(V)\times W$ with orthogonal complement of dimension
$n-1$. So $w_n(-(\Rr\oplus \lambda ))$, the Stiefel-Whitney class of
this orthogonal complement, is zero.
This argument, going back to \cite[Proposition 2.1]{CH},
could be used to replace the reference to
Lemma \ref{umkehr} in the proof above.
However, the argument involving the Euler class has the advantage
that it provides some information on the solution set, as is
explained in Appendix B.
}
\Rem{Taking a dual viewpoint, we may lift from $Y$ to the
projective bundle $P(\Rr\oplus\lambda )$ of the vector bundle
$\Rr\oplus\lambda$ and, by restricting to the fibrewise Hopf bundle
$H_{\Rr\oplus\lambda}\subseteq \Rr\oplus\lambda$, get
a section $s^*$ of $H_{\Rr\oplus\lambda}\otimes W$.

The direct image homomorphism
$$
H^{n+1}(P(\Rr\oplus \lambda );\,\Ff_2) \to  H^n(Y;\,\Ff_2)
$$
maps $e(H_{\Rr\oplus\lambda}\otimes W)$ to $w_n(-(\Rr\oplus\lambda ))$.
(See Lemma \ref{dual}.)
A zero of $s^*$ gives a point $x\in S(V)$ and 
$(\alpha ,\beta )\in \Rr^2-\{ (0,0)\}$ such that
$\alpha f(x)+\beta f(-x)=0$, or, in other words,
$f(x)$ and $f(-x)$ are parallel.
}
\Sect{Parallel lines}
The case in which $m=0$ and $f$ is a map $\Rr \to \Rr^2$ is special
and elementary. 
\Prop{Suppose that $f: \Rr \to\Rr^2$ is a continuous map.
Then there exist $4$ points $x_0<y_0<y_1<x_1$ in $\Rr$ such that
$f(x_1)-f(x_0)$ and $f(y_1)-f(y_0)$ are parallel.
}
\begin{proof}
If the image of $f$ is contained in some (affine) line, we can
choose the $4$ points arbitrarily.
If not, choose $3$ points $x_0 < z < x_1$ in $\Rr$ such that
$f(x_0),\, f(z),\, f(x_1)$ are not collinear. Apply the Intermediate
Value theorem to the function $\phi : [x_0,x_1]\to\Rr$ given by
the inner product $\rho (y) =\langle f(y)-f(x_0),e\rangle$ for some 
unit vector $e\in\Rr^2$ perpendicular to $f(x_1)-f(x_0)$.
Since $\rho (x_0)=0=\phi (x_1)$ and $\rho (z)\not=0$, there are points
$y_0$, $y_1$ such that $x_0<y_0<z<y_1<x_1$ such that $\rho (y_0)=
\rho (y_1)$.
\end{proof}
Assume now that $m\geq 1$. We consider as in Section \ref{prelude}
a continuous map $f: V=\Rr^{m+1}\to W=\Rr^{n+1}$.

The manifold $\hat X=S(V)\times S(V\oplus V)$, of dimension $3m+1$,
has a free action of the group $G=\{ \pm 1\} \times \Zz /2\Zz$ 
generated by the commuting involutions $(-1,0):
(x,(u,v))\mapsto (x,(-u,-v))$
and $(1,1):(x,(u,v))\mapsto (-x,(v,u))$.

The orbit space $\hat Y=\hat X/G$ can be identified with the real 
projective bundle $P(V\oplus (H\otimes V))$ by mapping the orbit
of $(x,(u,v))$ to $[u+v, x\otimes (u-v)]$ in the fibre over 
$[x]\in P(V)$.

There are two important line bundles over $\hat Y$:
the pullback $\lambda$ of the Hopf bundle $H$ over $P(V)$
and the Hopf bundle $\mu$ of the projective bundle.
The cohomology ring of $\hat Y$ is described in terms
of the mod $2$ Euler classes $t=e(\lambda )$ and $x=e(\mu )$
as
$$
H^*(\hat Y;\, \Ff_2)=\Ff_2[t,x]/(t^{m+1}, (x(t+x))^{m+1})
$$
(because $e(\mu\otimes (\lambda\otimes V))=
e(\mu\otimes\lambda )^{m+1}$).
Put $y=x(t+x)=e((\lambda\otimes\mu) \oplus \mu )$.
As $\Ff_2$-vector space, $H^*(\hat Y;\, \Ff_2)$ has 
a basis $t^iy^j,\, t^iy^jx$, $0\leq i,j\leq m$.

The bundles $\lambda$ and $\mu$ are associated with the
representations $G\to \O (\Rr )$:
$$
(a,c)\mapsto (-1)^c,\quad (a,c)\mapsto a,
$$
respectively.
(Isomorphisms are given by
$\hat X\times_G \Rr \to \lambda$:
$(x,(u,v),t)\mapsto tx$ and
$\hat X \times_G \Rr \to \mu$:
$((x, (u,v)), t)\mapsto t(u+v, x\otimes (u-v))$.)

Choose $\delta$, $0<\delta < \frac{1}{2}$.
There is a basic map (introduced, implicitly, in \cite{FH})
$$ 
\hat\phi :
\hat X= S(V)\times S(V\oplus V) \to V\oplus V\oplus V\oplus V\oplus V,
$$
$$
(x,(u,v)) \mapsto 
(x+\delta u,x-\delta u,-x+\delta v, -x-\delta v).
$$

\begin{proof}[Proof of Theorem \ref{parallel}]
The $4$-tuple $\hat \phi (x,(u,v))$ determines a pair of $2$-element
subsets of $V$:
$$
\{ x+\delta u, -x+\delta v\},\,
\{ x-\delta u, -x-\delta v\},
$$
with one point in common if $u=0$ or $v=0$, and otherwise disjoint.

Now $f$ determines sections of $\lambda\otimes W$ and 
$\lambda\otimes\mu\otimes W$ by the equivariant maps 
$$
(x,(u,v))\mapsto  
f(x+\delta u)-f(-x+\delta v)+f(x-\delta u)-f(-x-\delta v)
$$
and
$$
(x,(u,v))\mapsto  
f(x+\delta u)-f(-x+\delta v)-f(x-\delta u)+f(-x-\delta v).
$$
Taking their sum and projecting from $W$ to the quotient
$W/H$,
we get a section $s$ of $(\lambda\otimes (\Rr\oplus\mu ))
\otimes (W/H )$ over $\hat Y\times P(W)$.

At a zero $([x,(u,v)],L)\in \hat Y\times P(W)$ of $s$,
the vectors
$f(x+\delta u)-f(-x+\delta v)+f(x-\delta u)-f(-x-\delta v)$
and 
$f(x+\delta u)-f(-x+\delta v)-f(x-\delta u)+f(-x-\delta v)$
lie in the $1$-dimensional subspace $L$ of $W$.
So their sum $2(f(x+\delta u)-f(-x+\delta v))$ and
difference $2(f(x-\delta u)-f(-x-\delta v))$ lie in $L$.
Thus the vectors $f(x+\delta u)-f(-x+\delta v)$ and
$f(x-\delta u)-f(-x-\delta v)$ are parallel.

Now the image of the Euler class 
$$
e(\lambda\otimes (\Rr\oplus\mu )
\otimes (W/H ) )\in H^{2n}(\hat Y\times P(W);\,\Ff_2)
$$
in $H^n(\hat Y;\, \Ff_2)$ under the direct image homomorphism is 
$w_n(-(\lambda\otimes (\Rr\oplus\mu )))$.

To compute this Stiefel-Whitney class, note that
$$
1+t+y =w((\lambda\otimes\mu )\oplus \mu )=w(\lambda\otimes\mu )w(\mu)
=(1+t+x)(1+x).
$$
So
$$
w(-(\lambda\otimes (\Rr\oplus\mu ))=
(1+t)^{-1}(1+t+x)^{-1}=(1+t)^{-1}(1+t+y)^{-1}(1+x)
$$
$$
(1+t)^{-2}(1+y(1+t)^{-1})^{-1}(1+x)
=\sum_{j\geq 0} (1+t)^{-(j+2)}y^j(1+x)
$$
$$
=\sum_{0\leq i,j\leq m} \binom{-j-2}{i} t^iy^j(1+x)
=\sum_{i,j\leq m}\binom{2^r-j-2}{i}t^iy^j(1+x).
$$
The coefficient of $t^{2^r-m-2}y^mx$ is $1$,
and $2^r-m-2\leq m$.

With $n=m+2^r-1$, we have shown that the Euler class is
non-zero, so that the section $s$ has a zero and the
assertion (b) has been proved.

We also have a section $z$ of $\lambda$ determined by
the equivariant map
$$
(x,(u,v)) \mapsto \| u\|^2-\| v\|^2.
$$
At a zero of $z$, $u$ and $v$ are both non-zero.

The section $z\oplus s$ of 
$\lambda \oplus ((\mu\oplus (\lambda\otimes \mu ))\otimes (W/H))$
will (by Lemma \ref{umkehr}) have a zero if
$$
e(\lambda )w_n(-\mu -(\lambda\otimes\mu ))
=tw_n(-\mu -(\lambda\otimes\mu ))\in H^{n+1}(\hat Y;\,\Ff_2)
$$
is non-zero.
From the calculation above this is the case if $t(t^{2^r-m-2}y^mx)
\not=0$, which is true provided that $2^r-m-1\leq m$, that is,
$m+1\not=2^{r-1}$. This completes the proof of assertion (a).
\end{proof}
We discuss next the relation between these techniques and the
method used in \cite{FH}.
Writing $\O (V)$ for the orthogonal group of $V$,
let $G=\{ \pm 1\}\times \Zz /2\Zz$ act on 
$\hat X_0=\O (V)\times S(V\oplus V)$ by the involutions
$(g, (u,v))\mapsto (g, (-u,-v))$ and 
$(g, (u,v))\mapsto (-g, (u,v))$.
The orbit space $\hat Y_0=\hat X_0/G$ is the product
$\PO (V)\times P(V\oplus V)$ of the projective
orthogonal group of $V$ and the projective space.

Choosing a vector $e\in S(V)$ we have a $G$-map
$$
\hat X_0=\O (V)\times S(V\oplus V)\to \hat X =S(V)\times S(V\oplus V):
$$
$$
(g, (u,v)) \mapsto (ge, \frac{1}{\sqrt{2}}(u+gv,u-gv)).
$$
Under the quotient map $\hat Y_0 \to \hat Y$,
the line bundles $\lambda$ \and $\mu$ pullback to the
line bundle $\lambda_0$ associated with the double cover
$\O (V)\to\PO (V)$ and the Hopf line bundle $\mu_0$ over
$P(V\oplus V)$.
\Prop{Let $q\geq 0$ be the largest integer such that
$2^q$ divides $m+1$. 
Then
$w_n(-(\lambda_0\otimes (\Rr\oplus\mu_0))\not=0$ for
$n\leq 2m+2^q$.
}
\begin{proof}
It is classical that the Euler class $t_0=e(\lambda_0)$
satisfies $t_0^{2^q-1}\not=0$, $t_0^{2^q}=0$;
see, for example, \cite[Proposition A.1]{ivt}.
And $x_0=e(\mu_0)$ satisfies $x_0^{2m+1}\not=0$,
$x_0^{2m+2}=0$.
Thus $t_0$ and $t_0+x_0$ generate a subalgebra
of $H^*(\hat Y_0;\,\Ff_2)$ with the only relations
$t_0^{2^q}=0$ and $(t_0+x_0)^{2m+2}=0$.
(Notice that $(t_0+x_0)^{m+1}=(t_0^{t^q}+x_0^{2^q})^{(m+1)/2^q}
=x_0^{m+1}$.)
Now
$$
w(-(\lambda_0\otimes (\Rr\oplus \mu ))=
(1+t_0)^{-1}(1+t_0+x_0)^{-1}
=\sum_{i,\, j\geq 0} t_0^i(t_0+x_0)^j.
$$
So the assertion follows, because $2m+2^q=(2^q-1)+(2m+1)$.
\end{proof}
\Rem{The argument in \cite[Theorems 1 and 3]{FH} proceeds,
in effect, by using,
instead of the fact that $t_0^{2^q-1}\not=0$, a result of
the form $t_0^{2^{\alpha (q)}-1}\not=0$, 
where $\alpha (q)\leq q$ is expressed in terms of 
Hurwitz-Radon numbers.
For $q\geq 1$, $\alpha (q)\geq 1$ is chosen so that one has a
Clifford module multiplication $\Rr^{2^{\alpha (q)-1}+1}
\otimes V\to V$, and this determines
an embedding $P(\Rr^{2^{\alpha (q)-1}+1})
\into \PO (V)$, under which $\lambda_0$ restricts to the
Hopf line bundle. Hence $t_0^{2^{\alpha (q)-1}}\not=0$.
Now we recall a general $H$-space argument.
Consider the group multiplication
$\mu : \PO (V)\times\PO (V)\to \PO (V)$.
Since $\mu^*(t_0)=t_0\otimes 1+1\otimes t_0$,
$$
\mu^*(t_0^{2^{\alpha (q)}-1})
=(t_0\otimes 1+1\otimes t_0)^{2^{\alpha (q)}-1}
=\sum_{i+j=2^{\alpha (q)}-1} t_0^i\otimes t_0^j\, ,
$$
which is non-zero, because it contains the term
$t_0^{2^{\alpha (q)-1}-1}\otimes t_0^{2^{\alpha (q)-1}}$,
and hence $t_0^{2^{\alpha (q)}-1}\not=0$.
}
\Sect{Regular maps}
Let $X$ be the codimension one submanifold
$$
X=S(V)\times (S(V)\times S(V)) \subseteq \hat X
$$
of $\hat X$.
The quotient of $X$ by the group $G=\{ \pm 1\}\times\Zz /2\Zz$
is a submanifold $Y$ of $\hat Y$.
The quotient of $X$ by the free action of the group 
$K=(\{\pm 1\}\times \{ \pm 1\})\semidirect \Zz /2\Zz$ 
generated by the involutions
$((-1,1),0)$, $((1,-1),0)$ and $((1,1),1)$:
$(x,(u,v))\mapsto (x,(-u,v)),\, (x,(u,-v)),\, (-x, (v,u))$
is denoted by $Z$. It fibres over the real projective
space $P(V)$ on $V$ by the map induced by the projection
$(x,(u,v))\mapsto x$.
We have covers $X\to Y\to Z$.

There are various natural real vector bundles over $Z$:
the pull-back $\lambda$ of the Hopf line bundle
over $P(V)$, the line bundle $\alpha$ with fibre
$\Rr u\otimes \Rr v$ over the point determined
by $(x,(u,v))\in X$, and the $2$-dimensional vector
bundle $\zeta$ with fibre $\Rr u\oplus\Rr v$.
They are associated with the representations:
$$
(a,b; c)\mapsto (-1)^c;
$$
$$
(a,b; c)\mapsto ab;
$$
$$
(a,b; 0)\mapsto \left[\begin{matrix} a&0\\ 0&b\end{matrix}\right],
(a,b; 1)\mapsto \left[\begin{matrix} 0&a\\ b&0\end{matrix}\right],
$$
The determinant bundle $\Lambda^2\zeta$ is isomorphic to 
$\lambda\otimes\alpha$.

The mod $2$ Euler classes $t=e(\lambda )$,
$a=e(\alpha )$ and $y=e(\zeta )$ generate the
cohomology ring
$$
H^*(Z;\, \Ff_2) =\Ff_2[t,y,a]/(t^{m+1},\, y^{m+1},\, ta,\, 
y^ma,\, \ldots ,\, y^{m+1-i}a^i,\, \ldots ,\, a^{m+1})\, .
$$

The pullback of $\alpha$ to $Y$ is trivial. Indeed, $Y=S(\alpha )$.
The pullback of $\zeta$ to $Y$ splits as 
$(\lambda\otimes\mu)\oplus\mu$.
The kernel of the homomorphism
$$
H^*(Z;\,\Ff_2)\to H^*(Y;\,\Ff_2)
$$
is the ideal generated by $a$.

Let $\phi : X \to V\oplus V\oplus V\oplus V$ be the restriction
of $\hat\phi$.
Then $\phi (x,(u,v))$ determines $4$ distinct points 
$x\pm \delta u,\, -x\pm \delta v$ of $V$,
and a pair of disjoint $2$-element subsets
$$
\{ x+\delta u, x-\delta u\},\, \{ -x+\delta v,-x-\delta v\}
$$
of $V$.
\begin{proof}[Second proof of Theorem \ref{parallel}(a)] 
The map $f$ determines a section $s$ of $\zeta\otimes (W/H)$ 
over $Z\times P(W)$
by the equivariant map
$$
(x,(u,v)) \mapsto 
(f(x+\delta u)-f(x-\delta u)+f(-x+\delta v)-f(-x-\delta v),
$$
$$
f(x+\delta u)-f(x-\delta u)-f(-x+\delta v)+f(-x-\delta v)).
$$

The Euler class will
be non-zero if $w_n(-\zeta )\in H^n(Z;\,\Ff_2)$ is non-zero. 
We lift to $H^n(Y;\,\Ff_2)$ and compute
$$
(1+t+y)^{-1}=(1+t)^{-1}(1+y(1+t)^{-1})^{-1}
=\sum_{j\geq 0} (1+t)^{-(j+1)}y^j
$$
$$
=\sum_{0\leq i,j\leq m} \binom{-j-1}{i} t^iy^j
=\sum_{0\leq i,j\leq m} \binom{2^r-j-1}{i} t^iy^j\, .
$$
The coefficient of $t^{2^r-m-1}y^m$ is $1$,
and $2^r-m-1\leq m$, because $m+1\not=2^{r-1}$.
\end{proof}
\begin{proof}[Second proof of Corollary \ref{cor}]
We shall, in effect, reproduce the proof in \cite{mod2}
for this special case.

The equivariant map
$$
(x,(u,v)) \mapsto 
f(x+\delta u)-f(-x+\delta v)+f(x-\delta u)-f(-x-\delta v),
$$
determines a section of $\lambda\otimes W$, and, combined with the
section of $\zeta\otimes W$ already constructed,
a section of 
the $3(n-1)$-dimensional vector bundle
$(\lambda\oplus\zeta )\otimes (W/\gamma)$
over $Z\times G_2(W)$, where $\gamma$ is the canonical 
$2$-plane bundle
over the Grassman manifold $G_2(W)$ of $2$-dimensional subspaces 
of $W$.
The Euler class will be non-zero if
$w_{n-1}(-(\lambda\oplus\zeta ))\not=0$ (by \cite[Lemma 6.4]{tver}
or \cite[Proposition 2.1]{CH}).
Passing to the quotient of $H^*(Z;\,\Ff_2)$ by the ideal generated
by $a$, 
we look at the degree $n-1$ term in
$$
(1+t)^{-1}(1+t+y)^{-1}=\sum_{j\geq 0} (1+t)^{-(j+2)}y^j
=\sum_{0\leq i,j\leq m} \binom{2^r-j-2}{i} t^iy^j\, .
$$
The coefficient of 
$t^{2^r-m-2}y^m$ is equal to $1$.
Take $n-1=2^r-m-2+2m$.
\end{proof}
There is an equivalent spherical version of Corollary \ref{cor}
(which, in fact, is the form given in \cite[Theorem 6.16]{mod2}).
\Prop{\label{sphere}
Let $\tilde f: \Rr^{m+1}\to S(\Rr^{m+2^r+1})\subseteq \Rr^{m+2^r+1}$ 
be a continuous map, where $2^{r-1}\leq m+1<2^r$.

Then there exist $4$ distinct points $x_0,\, x_1,\, x_2,\, x_3\in
\Rr^{m+1}$ such that the vectors $\tilde f(x_0)$, $\tilde f(x_1)$, 
$\tilde f(x_2)$,
$\tilde f(x_3)$ are linearly dependent.
}
We give two proofs, one using the equivalence and one by a 
direct argument. Consider a general continuous map
$\tilde f : V \to S(\tilde W)$ to the unit sphere in
a Euclidean vector space $\tilde W$ of dimension $(n+1)+1$. 
\begin{proof}[First proof]
Let $i : W \to S(\Rr\oplus W)$ be the inclusion
$w\mapsto (1,w)/(1+\| w\|^2)^{1/2}$ of $W$ as the open
northern hemisphere of the sphere.
It has the property that $k+1$ vectors $w_0,\ldots ,w_k$ in
$W$ are affinely independent, that is, span a $k$-dimensional
affine subspace, if and only if $i(w_0),\ldots , i(w_k)$
are linearly independent in $\Rr\oplus W$, that is, span
a vector subspace which intersects
$S(\Rr\oplus W)$ in a $k$-dimensional sphere with centre
$0$. 

Now $\tilde f$ maps a sufficiently small
open neighbourhood of $0\in W$ into an open hemisphere in $S(\Rr\oplus W)$. So, by restricting $\tilde f$ to
a neighbourhood homeomorphic to $W$ and choosing an appropriate
splitting of $\tilde W=\Rr\oplus W$, 
we may assume that $\tilde f$ has the form $i\comp f$
for a continuous map $f : V\to W$. So the result
follows from Corollary \ref{cor}.
\end{proof}
\begin{proof}[Second proof]
The constructions in the second proof of Proposition \ref{cor}
applied to $\tilde f$ instead of $f$ give sections
of $\zeta\otimes\tilde W$ and $\lambda\otimes\tilde W$
and the map
$$
(x,(u,v)) \mapsto 
f(x+\delta u)+f(-x+\delta v)+f(x-\delta u)+f(-x-\delta v)
$$
gives a section of the trivial bundle with fibre $\tilde W$.
Combining the three maps we get a section of
$(\Rr\oplus\zeta\oplus\lambda )\otimes (\tilde W/\gamma )$
over $Z\times G_3(\tilde W)$, where $\gamma$ is now the
canonical $3$-dimensional bundle over the Grassmannian.
Its Euler class will be non-zero if
$w_{n+1-2}(-(\Rr\oplus\zeta\oplus\lambda ))\not=0$.
This is the same condition as in the affine case.
\end{proof}

\begin{appendix}
\Sect{}
We give an elementary proof of a special case of
\cite[Lemma 6.4]{tver}.
\Lem{\label{umkehr}
Let $\xi$ be a $2$-dimensional real vector bundle 
over a compact ENR $B$. Write $W=\Rr^{n+1}$.
Then the image of the Euler class $e(\xi\otimes (W/H))$
under the direct image homomorphism
$$
\pi_! : H^{2n}(B\times P(W);\,\Ff_2)\to H^{n}(B;\,\Ff_2)
$$
is equal to $w_n(-\xi )$. 
}
\begin{proof}
By universality, it suffices to verify the formula
when $\xi =\lambda_1\oplus \lambda_2$ is a
sum of line bundles. Write $t_i=e(\lambda_i)\in H^1(B;\,\Ff_2)$
and $x=e(H)\in H^1(P(W);\,\Ff_2)$.
Thus $w_{n-j}(W/H)=x^{n-j}$ and
$$
e((\lambda_1\oplus \lambda_2)\otimes (W/H))
=\prod_{i=1}^2 e(\lambda_i\otimes (W/H))
=\prod_{i=1}^2 (x^n+\ldots +t_i^jx^{n-j}+\ldots + t_i^n)
$$
in $H^*(B;\,\Ff_2)[x]/(x^{n+1})$.
We have to determine the coefficient of $x^n$
in this truncated polynomial ring,
or, equivalently, the coefficient of $T^{2n}$
in the formal power series
$$
\prod_{i=1}^2 (1+\ldots +t_i^jT^j+\ldots + t_i^nT^n)
=\prod_{i=1}^2 \frac{1-t_i^{n+1}T^{n+1}}{1-t_iT}
\in H^*(B;\,\Ff_2)[[T]].
$$
Now 
$$
\prod (1-t_iT)^{-1} = 
(1+w_1(\xi )T+ w_{2}(\xi )T^2)^{-1}
= 1+w_1(-\xi )T +\ldots + w_i(-\xi )T^i+ \ldots.
$$
So the coefficient of $T^{n}$ is $w_n(-\xi )$.
\end{proof}
And here is the dual version, a special case of 
the discussion preceding Proposition 4.2 in
\cite{borsuk}.
\Lem{\label{dual}
Let $\xi$ be a $2$-dimensional real vector bundle 
over a compact ENR $B$. Write $W=\Rr^{n+1}$.
Consider the projective bundle $\pi :P(\xi )\to B$.
Then the image of the Euler class $e(H\otimes W)$
under the direct image homomorphism
$$
\pi_! : H^{n+1}(P(\xi );\,\Ff_2)\to H^{n}(B;\,\Ff_2)
$$
is equal to $w_n(-\xi )$. 
}
\begin{proof}
The cohomology ring of the projective bundle is
$$
H^*(P(\xi );\,\Ff_2) =H^*(B)[t]/(t^2+w_1(\xi )t+w_2(\xi )),
$$
where $t=e(H)$, so that $e(H\otimes W)=t^{n+1}$. 
Induction on $n$, 
using the identity $w_{n+1}(-\xi )+w_1(\xi )w_n(-\xi )
+w_2(\xi )w_{n-1}(-\xi )=0$ for $n\geq 1$, shows
that $t^{n+1}=w_n(-\xi )t+w_2(\xi )w_{n-1}(-\xi )$.
The Umkehr $\pi_!$ picks out the coefficient of $t$.
\end{proof}
\Sect{}
Consider a continuous map $f: V=\Rr^{m+1}\to W=\Rr^{n+1}$
and write $\Sing (f)$ for the subspace
of $V^4$ consisting of those $4$-tuples $(x_0,x_1,x_2,x_3)$
such that $x_i\not=x_j$ for $i\not=j$
and the points $f(x_0),\, f(x_1),\, f(x_2),\, f(x_3)$ lie
in some $2$-dimensional affine subspace of $W$.
We shall establish a lower bound for the covering dimension
of this {\it singularity set} $\Sing (f)$.
\Prop{Suppose that $n+1=m+2^r$ where $2^{r-1}\leq m+1<2^r$.
Then $\Sing (f)$ contains a compact subspace with
covering dimension greater than or equal to $4(m+1)-(n-1)$.
}
\begin{proof}
We follow the notation of Section 4 where we used the inclusion
$(x,(u,v))\mapsto (x+\delta u,x-\delta u,-x+\delta v,-x-\delta v)$:
$$
X=S(V)\times (S(V)\times S(V)) \into V^4.
$$
Choose a $K$-equivariant closed tubular neighbourhood $\tilde X$
of $X$ in $V^4$, so small that, 
for all $(x_0,x_1,x_2,x_3)\in\tilde X$, we have
$x_i\not=x_j$ for $i\not=j$.
Then $\tilde Z=\tilde X/K$ is a compact manifold with boundary
and the inclusion $Z\into\tilde Z$ is a homotopy
equivalence. Let $\tilde\lambda$ and $\tilde\zeta$ denote
the natural extensions of $\lambda$ and $\zeta$ to $\tilde Z$.
The construction in the second proof of Corollary \ref{cor} 
in Section 4 gives a section, $s$ say, of the vector bundle
$(\tilde\lambda\oplus\tilde\zeta )\otimes (W/\gamma )$ over 
$\tilde Z\times G_2(W)$.
And the zero-set $\Zero (s)$ projects under
$\pi : \tilde Z\times G_2(W)\to \tilde Z$ to
the image of the compact space $\Sing (f)\cap \tilde X$ 
under the $8$-fold covering map $\tilde X\to\tilde Z$.

We know that the Stiefel-Whitney class 
$w_{n-1}(-(\tilde\lambda\oplus\tilde\zeta ))$ is non-zero.
By Poincar\'e duality for the manifold $\tilde Z$ of 
dimension $4(m+1)$, there is a class 
$$
a\in H^{4(m+1)-(n-1)}(\tilde Z,\partial\tilde Z;\,\Ff_2)
$$ 
such that
$$
a\cdot w_{n-1}(-(\tilde\lambda\oplus\tilde\zeta))=
\pi_!((\pi^*a)\cdot 
e((\tilde\lambda \oplus\tilde\zeta )\otimes (W/\gamma )))
$$
is non-zero.
Thus $(\pi^*a )\cdot
 e((\tilde\lambda \oplus\tilde\zeta )\otimes (W/\gamma ))$
is non-zero. It follows by standard arguments,
as, for example, in \cite[Proposition 2.7]{borsuk},
that $\pi^*a$ restricts to a non-zero class in the (representable)
cohomology of
$\Zero (s)$. Hence $a$ must restrict to a non-zero class in
the cohomology of the image $S$ of $\Sing (f)\cap \tilde X$
in $\tilde Z$. Therefore
the covering dimension of the compact space
$S$ is at least $4(m+1)-(n-1)$, and so, too,
is the covering dimension
of its $8$-fold cover $\Sing (f)\cap \tilde X$.
\end{proof}
\end{appendix}

\end{document}